\documentclass[12pt]{amsart}
\usepackage{amsmath, amsthm, amscd, amsfonts}

\newtheorem{theorem}{Theorem}[section]

\newtheorem{proposition}[theorem]{Proposition}
\newtheorem{corollary}[theorem]{Corollary}
\theoremstyle{definition}

\theoremstyle{remark}
\newtheorem{remark}[theorem]{Remark}
\numberwithin{equation}{section}

\begin{document}
\title{Hyers--Ulam--Rassias Stability of Generalized Derivations }
\author{Mohammad Sal Moslehian}
\address{Department of Mathematics, Ferdowsi University, P. O. Box 1159, Mashhad 91775,
Iran \newline and \newline Banach Mathematical Reseach Group
(BMRG)} \email{moslehian@ferdowsi.um.ac.ir}
\subjclass[2000]{Primary 39B82; Secondary 46H25, 39B52, 47B47.}
\keywords{Hyers--Ulam--Rassias stability, normed algebra, Banach
module, generalized derivation, $C^*$-algebra.}
\begin{abstract}
The generalized Hyers--Ulam--Rassias stability of generalized
derivations on unital Banach algebras into Banach bimodules is
established.
\end{abstract}
\maketitle
\section {Introduction}

One of interesting questions in the theory of functional
equations is the problem of the stability of functional equations
as follows: ``When is it true that a mapping satisfying a
functional equation approximately must be close to an exact
solution of the given functional equation?''

The first stability problem was raised by S. M. Ulam during his
talk at University of Wisconsin in 1940 \cite{ULA}:

{\small Given a group $G_1$, a metric group $(G_2, d)$ and a
positive number $\varepsilon$, does there exist a $\delta>0$ such
that if a mapping $f : G_1 \to G_2$ satisfies the inequality
$d(f(xy), f(x)f(y))<\delta$ for all $x, y \in G_1$ then there
exists a homomorphism $T : G_1 \to G_2$ such that $d(f(x),
T(x))<\varepsilon$ for all $x\in G_1$.}

Ulam's problem was partially solved by D. H. Hyers in 1941 in the
context of Banach spaces with $\delta=\varepsilon$ as the
following \cite{HYE}:

{\small Suppose that $E_1, E_2$ are Banach spaces and $f:E_1 \to
E_2$ is a mapping for which there exists $\varepsilon>0$ such that
$\|f(x+y)-f(x)-f(y)\|<\varepsilon$ for all $x, y \in E_1$. Then
there is a unique additive mapping $T:E_1\to E_2$ defined by
$Tx=\displaystyle{\lim_{n\to\infty}}\frac{f(2^nx)}{2^n}$ such
that $\|f(x)-T(x)\|<\varepsilon$ for all $x\in E_1$.}

Now assume that $E \sb 1$ and $E\sb 2$ are real normed spaces with
$E\sb 2$ complete, $f: E\sb 1\to E\sb 2$ is a mapping such that
for each fixed $x\in E\sb 1$ the mapping $t\mapsto f(tx)$ is
continuous on ${\mathbb R}$, and let there exist $\varepsilon\ge
0$ and $p\neq 1$ such that
$$\|f(x+y)-f(x)-f(y)\|\le \varepsilon(\|x\|\sp p+\|y\|\sp p)$$ for all $x, y \in E\sb 1$.

It was shown by Th. M. Rassias \cite{RAS1} for $p\in [0, 1)$ (and
indeed $p<1$) and Z. Gajda \cite{GAJ} following the same approach
as in \cite{RAS1} for $p>1$ that there exists a unique linear map
$T: E\sb 1 \to E\sb 2$ such that
$$\|f(x)-T(x)\|\le\frac{2\varepsilon}{|2^p-2|}\|x\|\sp p$$ for all
$x \in E\sb 1$. This phenomenon is called {\it
Hyers--Ulam--Rassias stability}. It is shown that there is no
analogue of Th. M. Rassias result for $p=1$ (see \cite{ GAJ,
R-S}).

In 1992, a generalization of the Rassias theorem was obtained by
G\u avruta as follows \cite{GAV}:

{\small Suppose $(G,+)$ is an abelian group, $E$ is a Banach
space and the so-called admissible control function
$\varphi:G\times G\to [0, \infty)$ satisfies
$$\widetilde{\varphi}(x, y):
=\frac{1}{2}\displaystyle{\sum_{n=0}^\infty}2^{-n}\varphi(2^n
x,2^n y)<\infty$$ for all $x,y\in G$. If $f : G \to E$ is a
mapping with $$\|f(x+y)-f(x)-f(y)\|\leq \varphi(x,y)$$ for all
$x, y \in G$, then there exists a unique mapping $T : G \to E$
such that $T(x+y)=T(x)+T(y)$ and
$\|f(x)-T(x)\|\leq\widetilde{\varphi}(x, x)$ for all $x,y\in G$.}

Since then several stability problems of various functional
equations have been investigated by many mathematicians. The
reader is referred to \cite{CZE, RAS2} for a comprehensive account
of the subject.

Generalized derivations were first appeared in the context of
operator algebras \cite{MAT}. Later, these were introduced in the
framework of pure algebra \cite{HVA}. There is also other
generalization of the notion of derivation which is called
$(\sigma,\tau)$-derivation; cf. \cite{M-M}.

Let ${\mathcal A}$ be an algebra and let ${\mathcal X}$ be an
${\mathcal A}$-bimodule. A linear mapping $\mu:{\mathcal A}\to
{\mathcal X}$ is called a {\it generalized derivation} if there
exists a derivation (in the usual sense) $\delta:{\mathcal A}\to
{\mathcal X}$ such that $\mu(ab)=a\mu(b)+\delta(a)b$ for all
$a,b\in {\mathcal A}$. Familiar examples are the derivations from
${\mathcal A}$ to ${\mathcal X}$ and all so-called inner
generalized derivations i.e. those are defined by
$\mu_{x,y}(a)=xa-ay$ for fixed arbitrary elements $x,y\in
{\mathcal X}$. Moreover, every right multiplier (i.e., an additive
map $h$ of ${\mathcal A}$ satisfying $h(xy) = h(x)y$ for all $x,
y \in {\mathcal A}$) is a generalized derivation.

The stability of derivations was studied by C.-G. Park in
\cite{PAR1, PAR2}. A discussion of stability of the so-called
$(\sigma-\tau)$-derivations and a study of the so-called
generalized $(\theta,\phi)$-derivations are given in \cite{MOS1}
and \cite{B-M1}, respectively. The present paper is devoted to
study of stability of generalized derivations. The results of
this paper are a generalization of those of Park's papers
\cite{PAR1, PAR2}.

Throughout the paper $A$ denotes a unital normed algebra with
unit $1$ and ${\mathcal X}$ is a unit linked  Banach ${\mathcal
A}$-bimodule in the sense that $1x=x1=x$ for all $x\in {\mathcal
X}$.

\section{Main Results.}

Our aim is to establish the generalized Hyers--Ulam--Rassias
stability of generalized derivations. We extend main results of
C.-G. Park \cite{PAR1} to generalized derivations from a unital
normed algebra to a unit linked Banach ${\mathcal A}$-bimodule.
We apply the direct method which was first devised by D. H. Hyers
\cite{HYE} to construct an additive function from an approximate
one and use some ideas of \cite{MOS1} and \cite{PAR2}.

\begin{theorem} Suppose $f:{\mathcal A}\to {\mathcal X}$ is a mapping with $f(0)=0$ for which
there exist a map $g:{\mathcal A}\to {\mathcal X}$ and a function
$\varphi: {\mathcal A}\times {\mathcal A}\times {\mathcal A}\times
{\mathcal A}\to [0,\infty)$ such that
\begin{eqnarray}
\widetilde{\varphi}(a,b,c,d):=\frac{1}{2}\displaystyle{\sum_{n=0}^\infty}2^{-n}\varphi(2^na,2^nb,2^nc,2^nd)<\infty
\end{eqnarray}
\begin{eqnarray}
\|f(\lambda a+\lambda b+cd)-\lambda f(a)-\lambda f(b)-cf(d)-g(c)d\|\leq \varphi(a,b,c,d)
\end{eqnarray}
for all $\lambda\in {\mathbb T}=\{\lambda\in{\mathbb C}:
|\lambda|=1\}$ and all $a,b,c,d\in {\mathcal A}$. Then there
exists a unique generalized derivation $\mu:{\mathcal A}\to
{\mathcal X}$ such that
\begin{eqnarray}
\|f(a)-\mu(a)\|\leq\widetilde{\varphi}(a,a,0,0)
\end{eqnarray}
for all $a\in {\mathcal A}$.
\end{theorem}

\begin{proof} Setting $c=d=0$ and $\lambda=1$ in $(2)$, we have
\begin{eqnarray}
\|f(a+b)-f(a)-f(b)\|\leq\varphi(a,b,0,0),
\end{eqnarray}
for all $a,b\in {\mathcal A}$. Now we use the Th. M. Rassias
method on inequality $(2.4)$ (see \cite{GAV} and \cite{MOS2}).
One can use induction on $n$ to show that
\begin{eqnarray}
\|\frac{f(2^na)}{2^n}-f(a)\|\leq\frac{1}{2}\displaystyle{\sum_{k=0}^{n-1}}2^{-k}\varphi(2^ka,2^ka,0,0)
\end{eqnarray}
for all $n\in{\mathbb N}$ and all $a\in {\mathcal A}$, and that
\begin{eqnarray*}
\|\frac{f(2^na)}{2^n}-\frac{f(2^ma)}{2^m}\|\leq\frac{1}{2}\sum_{k=m}^{n-1}2^{-k}\varphi(2^ka,2^ka,0,0)
\end{eqnarray*}
for all $n>m$ and all $a\in {\mathcal A}$. It follows from the
convergence $(2.1)$ that the sequence $\{\frac{f(2^na)}{2^n}\}$
is Cauchy. Due to the completeness of ${\mathcal X}$, this
sequence is convergent. Set
\begin{eqnarray}
\mu(a) : = \displaystyle{\lim_{n\to\infty}}\frac{f(2^na)}{2^n}
\end{eqnarray}
Putting $c=d=0$ and replacing $a, b$ by $2^na, 2^nb$,
respectively, in $(2.2)$, we get
\begin{eqnarray*}
\|2^{-n}f(2^n(\lambda a+\lambda b))-2^{-n}\lambda
f(2^na)-2^{-n}\lambda f(2^nb)\|\leq 2^{-n}\varphi(2^na,2^nb,0,0).
\end{eqnarray*}
Taking the limit as $n\to\infty$ we obtain
\begin{eqnarray}
\mu(\lambda a+\lambda b)=\lambda \mu(a)+\lambda \mu(b),
\end{eqnarray}
for all $a,b\in {\mathcal A}$ and all $\lambda \in {\mathbb T}$.

Next, let $\gamma=\theta_1+{\bf i}\theta_2\in {\mathbb C}$ where
$\theta_1, \theta_2\in{\mathbb R}$. Let
$\gamma_1=\theta_1-[\theta_1], \gamma_2=\theta_2-[\theta_2]$.
Then $0\leq \gamma_i<1, ~(1\leq i \leq 2)$ and by using Remark
2.2.2 of \cite{MUR} one can represent $\gamma_i$ as
$\gamma_i=\frac{\lambda_{i,1} + \lambda_{i,2}}{2}$ in which
$\lambda_{i,j}\in {\mathbb T}, ~~(1\leq i,j \leq 2)$. Since $\mu$
satisfies $(2.7)$ we infer that
\begin{eqnarray*}
\mu(\gamma x)&=& \mu(\theta_1x)+{\bf i} \mu(\theta_2x)\\
&=&[\theta_1]\mu(x)+\mu(\gamma_1x)+{\bf i}\big( [\theta_2]\mu(x)+\mu(\gamma_2x)\big )\\
&=&\big( [\theta_1]\mu(x)+\frac{1}{2}\mu(\lambda_{1,1}x + \lambda_{1,2}x)\big )+{\bf i} \big ( [\theta_2]\mu(x)+\frac{1}{2}\mu(\lambda_{2,1}x + \lambda_{2,2}x)\big )\\
&=&\big([\theta_1]\mu(x)+\frac{1}{2}\lambda_{1,1}\mu(x) +\frac{1}{2}\lambda_{1,2}\mu(x)\big )+{\bf i}([\theta_2]\mu(x)+\frac{1}{2}\lambda_{2,1}\mu(x) + \frac{1}{2}\lambda_{2,2}\mu(x)\big )\\
&=& \theta_1\mu(x)+{\bf i}\theta_2\mu(x)\\
&=&\gamma \mu(x).
\end{eqnarray*}
for all $x\in {\mathcal A}$. So $\mu$ is ${\mathbb C}$-linear.

Moreover, it follows from $(2.5)$ and $(2.6)$ that
$\|f(a)-\mu(a)\|\leq \widetilde{\varphi}(a,a,0,0)$ for all $a\in
{\mathcal A}$. It is known that additive mapping $\mu$ satisfying
$(2.3)$ is unique \cite{B-M2}.

Putting $\lambda=0$, $x=y=0$ and replacing $c,d$ by $2^nc, 2^nd$,
respectively, in $(2.2)$ we obtain
$$\|f(2^{2n}cd)-2^ncf(2^nd)-2^ng(2^nc)d\|\leq\varphi(0,0,2^nc,2^nd),$$
whence
\begin{eqnarray}
\|2^{-2n}f(2^{2n}cd)-2^{-n}cf(2^nd)-2^{-n}g(2^nc)d\|\leq
2^{-2n}\varphi(0,0,2^nc,2^nd)
\end{eqnarray}
Put $d=1$ in $(2.8)$. By $(2.6)$,
$\displaystyle{\lim_{n\to\infty}}2^{-2n}f(2^{2n}a)=\mu(a)$ and by
the convergence of series $(2.1)$,
$\displaystyle{\lim_{n\to\infty}}2^{-2n}\varphi(0,0,2^nc,2^nd)=0$.
Hence the sequence $\{2^{-n}g(2^nc)\}$ is convergent. Set
$\delta(c) : = \displaystyle{\lim_{n\to\infty}}2^{-n}g(2^nc),~
c\in {\mathcal A}$. Let $n$ tend to $\infty$ in $(2.8)$,. Then
\begin{eqnarray}
\mu(cd)=c\mu(d)+\delta(c)d
\end{eqnarray}
Next we claim that $\delta$ is a derivation. Put $d=1$ in $(2.9)$.
Then $\delta(c)=\mu(c)-c\mu(1)$. Hence $\delta$ is linear.
Further,
\begin{eqnarray*}
\delta(c_1c_2) &=& \mu(c_1c_2)-c_1c_2\mu(1)\\
&=&(c_1\mu(c_2)+\delta(c_1)c_2)-c_1c_2\mu(1)\\
&=&c_1\mu(c_2)+(\mu(c_1)-c_1\mu(1))c_2-c_1c_2\mu(1)\\
&=&c_1(\mu(c_2)-c_2\mu(1))+(\mu(c_1)-c_1\mu(1))c_2\\
&=&c_1\delta(c_2)+\delta(c_1)c_2.
\end{eqnarray*}
Thus $\delta$ satisfies Leibnitz' rule. It then follows from
$(2.9)$ that $\mu$ is a generalized derivation.
\end{proof}

\begin{remark} The significance of functional equation $(2.2)$ is that the required derivation $\delta$ is naturally constructed.
In other words, we do not need any additional functional
inequality for existence of $\delta$.
\end{remark}

\begin{remark} Due to ${\mathcal A}$ is unital, the mapping $\delta$ appeared in the definition of
generalized derivation is unique. In fact,
$\delta(a)=\mu(a)-a\mu(1)$.
\end{remark}

\begin{corollary} Suppose that $f:{\mathcal A}\to {\mathcal X}$ is a mapping with $f(0)=0$ for which there exist constants $ \beta\geq 0$ and $p<1$ such that
\begin{eqnarray*}
\|f(\lambda a+\lambda b+cd)-\lambda f(a)-\lambda f(b)-cf(d)-g(c)d\|\leq\beta(\|a\|^p+\|b\|^p+\|c\|^p+\|d\|^p)
\end{eqnarray*}
for all $\lambda\in {\mathbb T}$ and all $a,b,c,d\in {\mathcal A}$.\\
Then there is a unique generalized derivation $\mu:{\mathcal
A}\to {\mathcal X}$ such that
\begin{eqnarray*}
\|f(a)-\mu(a)\|\leq\frac{\beta\|a\|^p }{1-2^{p-1}}
\end{eqnarray*}
for all $a\in {\mathcal A}$.
\end{corollary}
\begin{proof} Put $\varphi(a,b,c,d)=\beta(\|a\|^p+\|b\|^p+\|c\|^p+\|d\|^p)$ in
Theorem 2.1.
\end{proof}

\begin{proposition} Suppose that $f:{\mathcal A}\to {\mathcal X}$ is a mapping
with $f(0)=0$ for which there exists a function $\varphi:
{\mathcal A}\times {\mathcal A}\times {\mathcal A}\times
{\mathcal A}\to [0,\infty)$ such that
\begin{eqnarray*}
\widetilde{\varphi}(a,b,c,d):=\frac{1}{2}\displaystyle{\sum_{n=0}^\infty}2^{-n}\varphi(2^na,2^nb,2^nc,2^nd)<\infty
\end{eqnarray*}
\begin{eqnarray*}
\|f(\lambda a+\lambda b+cd)-\lambda f(a)-\lambda f(b)-cf(d)-g(c)d\|\leq \varphi(a,b,c,d)
\end{eqnarray*}
for $\lambda=1,{\bf i}$ and for all $a,b,c,d\in {\mathcal A}$. If
for each fixed $a\in {\mathcal A}$ the function $t\mapsto f(ta)$
is continuous on ${\mathbb R}$ then there exists a unique
generalized derivation $\mu:{\mathcal A}\to {\mathcal X}$ such
that $\|f(a)-\mu(a)\|\leq\widetilde{\varphi}(a,a,0,0)$ for all
$a\in {\mathcal A}$.
\end{proposition}
\begin{proof} Put $c=d=0$ and $\lambda=1$ in $(2.2)$. It follows from
the proof of Theorem 2.1 that there exists a unique additive
mapping $\mu:{\mathcal A}\to {\mathcal X}$ given by
$\mu(a)=\displaystyle{\lim_{n\to\infty}}\frac{f(2^na)}{2^n}, a\in
{\mathcal A}$. By the same reasoning as in the proof of the main
theorem of \cite{RAS1}, the mapping $\mu$ is ${\mathbb R}$-linear.

Assuming $b=c=d=0$ and $\lambda={\bf i}$, it follows from $(2.2)$
that $\|f({\bf i}a)-{\bf i}f(a)\|\leq\varphi(a,0,0,0), a\in
{\mathcal A}$. Hence $\frac{1}{2^n}\|f(2^n{\bf i}a)-{\bf
i}f(2^na)\|\leq\varphi(2^na,0,0,0)$ for all $n\in N$ and $a\in
{\mathcal A}$. The right hand side tends to zero as $n\to
\infty$, so that
\begin{eqnarray*}
\mu({\bf i}a)&=&\lim_{n\to\infty}\frac{f(2^n{\bf i}a)}{2^n}\\
&=&\lim_{n\to\infty}\frac{if(2^na)}{2^n}\\
&=&{\bf i}\mu(a)
\end{eqnarray*}
for all $a\in {\mathcal A}$. For each $\lambda\in{\mathbb C},
\lambda=r_1+{\bf i}r_2 ~~(r_1,r_2\in{\mathbb R})$. Hence
\begin{eqnarray*}
\mu(\lambda a)&=&\mu(r_1a+{\bf i}r_2a)=r_1\mu(a)+r_2\mu({\bf i}a)\\
&=&r_1\mu(a)+{\bf i}r_2\mu(a)=(r_1+{\bf i}r_2)\mu(a)\\
&=&\lambda\mu(a).
\end{eqnarray*}
Thus $\mu$ is ${\mathbb C}$-linear. That $\mu$ is a generalized
derivation can be deduced in the same fashion as in the proof of
Theorem 2.1.
\end{proof}

\begin{proposition} Let ${\mathcal A}$ be a unital $C^*$-algebra. Suppose that $f:{\mathcal A}\to {\mathcal X}$ is a
mapping with $f(0)=0$ for which there exists a function $\varphi:
{\mathcal A}\times {\mathcal A}\times {\mathcal A}\times
{\mathcal A}\to [0,\infty)$ such that
\begin{eqnarray*}
\widetilde{\varphi}(a,b,c,d):=\frac{1}{2}\displaystyle{\sum_{n=0}^\infty}2^{-n}\varphi(2^na,2^nb,2^nc,2^nd)<\infty
\end{eqnarray*}
\begin{eqnarray*}
\|f(\lambda a+\lambda b+cd)-\lambda f(a)-\lambda f(b)-cf(d)-g(c)d\|\leq \varphi(a,b,c,d)
\end{eqnarray*}
\begin{eqnarray}
\|f(2^nu^*)-f(2^nu)^*\|\leq \varphi(2^nu,2^nu,0,0)
\end{eqnarray}
for all $\lambda\in {\mathbb T}$, all $a,b,c,d\in {\mathcal A}$,
all nonnegative integers $n$ and all unitaries $u$ in ${\mathcal
A}$. Then there exists a unique generalized derivation
$\mu:{\mathcal A}\to {\mathcal X}$ such that
$\|f(a)-\mu(a)\|\leq\widetilde{\varphi}(a,a,0,0)$ for all $a\in
{\mathcal A}$.
\end{proposition}
\begin{proof} It follows from the proof of Theorem 2.1 that there exists a unique generalized derivation $\mu:{\mathcal A}\to {\mathcal X}$ given by $\mu(a)=\displaystyle{\lim_{n\to\infty}}\frac{f(2^na)}{2^n}, a\in {\mathcal A}$ satisfying $(2.3)$.\\
Using $(2.10)$, we have
\begin{eqnarray*}
\|2^{-n}f(2^nu^*)-2^{-n}f(2^nu)^*\|\leq 2^{-n}\varphi(2^nu,2^nu,0,0).
\end{eqnarray*}
Letting $n\to\infty$ we conclude that $\mu(u^*)=\mu(u)^*$. Since
$\mu$ is linear and every element of a $C^*$-algebra can be
represented as a linear combination of unitaries \cite{MUR}, we
deduce that $\mu(a^*)=\mu(a)^*$.
\end{proof}

Now let ${\mathcal A}$ be a unital Banach algebra. The mapping
$f:{\mathcal A}\to {\mathcal A}$ is called an {\it approximately
generalized derivation} if $f(0)=0$ and there exist a positive
number $\varepsilon$ and a mapping $g:{\mathcal A}\to {\mathcal
A}$ such that \[ \|f(\lambda a+\lambda b+cd)-\lambda f(a)-\lambda
f(b)-cf(d)-g(c)d\|\leq \varepsilon\] for all $\lambda\in {\mathbb
T}$ and all $a,b,c,d\in {\mathcal A}$.

\begin{theorem} Let ${\mathcal A}$ be a unital Banach algebra
and $f:{\mathcal A}\to {\mathcal A}$ be an approximately
generalized derivation with the corresponding mapping $g$. Then
$f$ is a generalized derivation and $g$ is a derivation.
\end{theorem}
\begin{proof}
Put $\varphi(a, b) = \varepsilon$ in Theorem 2.1. Then we get a
generalized derivation $\mu$ defined by $\mu(a):=
\lim_{n\to\infty}\frac{f(2^n a)}{2^n}$ such that
\[ \|\mu(a) - f(a)\| \leq \varepsilon\]
for all $a\in{\mathcal A}$. We have
\begin{eqnarray*}
\|2^n(f(2^ma) - 2^mf(a))\| &\leq& \|2^n1f(2^ma) - g(2^n1)2^ma -
f((2^n1)(2^ma))\|\\
&&+\|f((2^n1)(2^ma)) - g(2^n1)2^ma - 2^{n+m}1f(a)\|\\
&\leq& \varepsilon + \|f((2^n1)(2^ma)) - g(2^n1)2^ma -
2^{n+m}1f(a)\|\\
&\leq& \varepsilon + \|f((2^n1)(2^ma))- \mu((2^n1)(2^ma))\|\\
&&+\|\mu((2^n1)(2^ma)) - 2^{n+m}1f(a) - g(2^n1)2^ma\|\\
&\leq& 2\varepsilon + \|\mu((2^n1)(2^ma)) - 2^{n+m}1f(a) - g(2^n1)2^ma\|\\
&\leq& 2\varepsilon + 2^m \|\mu(2^n1a) - f(2^n1a)\|\\
&&+ 2^m \|f(2^n1a)
- 2^n1f(a) - g(2^n1)a\|\\
&\leq& (2 + 2^{m+1}) \varepsilon,
\end{eqnarray*}
for all nonnegative integers $m,n$ and all $a\in{\mathcal A}$.
Fix $m$ and let $n$ tend to $\infty$ in the following inequality
\begin{eqnarray*}
\|f(2^ma) - 2^mf(a)\|\leq \frac{2 + 2^{m+1}}{2^n} \varepsilon.
\end{eqnarray*}
Then $f(2^ma) = 2^mf(a)$ for all $m$ and all $a\in{\mathcal A}$.
Therefore $\mu(a) = \lim_{m\to\infty}\frac{f(2^m a)}{2^m} = f(a)$
for all $a\in{\mathcal A}$.
\end{proof}

\textbf{Acknowledgement.} The author would like to thank the
referees for their useful comments.

\end{document}